\documentclass{amsart}
\usepackage{graphicx}
\usepackage{amssymb}
\begin{document}

\centerline{\huge\bf Some Mathematical Problems Behind}

\medskip
\centerline{\huge\bf Lattice-Based Cryptography}

\bigskip
\centerline{\large Chuanming Zong}

\vspace{2cm}\noindent
In 1994, P. Shor discovered quantum algorithms which can break both the RSA cryptosystem and the ElGamal cryptosystem. In 2007, D-Wave demonstrated the first quantum computer. These events and further developments have brought a crisis to secret communication. In 2016, the National Institute of Standards and Technology (NIST) launched a global project to solicit and select a handful of encryption algorithms with the ability to resist quantum computer attacks. In 2022, it announced four candidates, CRYSTALS-Kyber, CRYSTALS-Dilithium, Falcon and Sphincs$+$ for post-quantum cryptography standards. The first three are based on lattice theory and the last on the Hash function. It is well known that the security of the lattice-based cryptosystems relies on the computational complexity of the shortest vector problem (SVP), the closest vector problem (CVP) and their generalizations. In fact, the SVP is a ball packing problem and the CVP is a ball covering problem. Furthermore, both SVP and CVP are equivalent to arithmetic problems for positive definite quadratic forms. Therefore, post-quantum cryptography provides unprecedented opportunities for mathematicians to make contributions in modern technology. This paper will briefly review the mathematical problems on which the lattice-based cryptography is built up, so that mathematicians can see that they are indeed in the game.

\bigskip
\noindent
{\large\bf 1. Introduction}

\smallskip\noindent
In 1976, W. Diffie and M. E. Hellman \cite{Diff76} proposed the principle of public key cryptography. One year later, the first public key cryptosystem RSA was invented by R. L. Rivest, A. Shamir and L. Adleman \cite{Rive78}.

These events not only inaugurated a new era in secret communication, but also marked the birth of mathematical cryptography. In the following years, several other mathematical cryptosystems have been discovered, including the discrete logarithm cryptosystem invented by T. ElGamal \cite{ElGa85} in 1985 and the elliptic curve cryptosystem ECC designed by V. S. Miller \cite{Mill85} and N. Koblitz \cite{Kobl87} in 1986, respectively.

In 1994, Shor \cite{Shor94} discovered polynomial time quantum algorithms to deal with the discrete logarithm problem and the factorization problem.
It follows that once there is a functioning quantum computer, Shor's algorithms could break both the RSA cryptosystem and the ElGamal cryptosystem.
At that time, quantum computer was still a fiction concept. People were not much bothered, excepting handful cryptographers.

From 1996 to 1998, three lattice-based cryptosystems were discovered, the Ajtai-Dwork by M. Ajtai and C. Dwork \cite{Ajta97}, the GGH by O. Goldreich, S. Goldwasser and S. Halevi \cite{Gold97}, and the NTRU by J. Hoffstein, J. Pipher and J. H. Silverman \cite{Hoff98}. Afterwards, several other cryptosystems have been created, including the LWE by O. Regev \cite{Rege05} in 2005 and the FHE by C. Gentry \cite{Gent09} in 2009. In the past half century, mathematical cryptography has played a crucial role in the modern technology of computers and the internet. It has been developed into an active interdisciplinary research field between mathematics and cryptography.

\smallskip
Lattice is a mathematical concept introduced by Gauss at the beginning of the 19th century and further developed by Minkowski and many others. Let ${\bf a}_1$, ${\bf a}_2$, $\ldots$, ${\bf a}_n$ be $n$ linearly independent vectors in the $n$-dimensional Euclidean space $\mathbb{E}^n$. We call
$$\Lambda =\left\{z_1{\bf a}_1+z_2{\bf a}_2+\ldots +z_n{\bf a}_n:\ z_i\in \mathbb{Z}\right\}$$
an $n$-dimensional lattice and call $\{ {\bf a}_1, {\bf a}_2, \ldots, {\bf a}_n\}$ a basis of the lattice $\Lambda$.

In 1998, the first quantum computer models appeared at Oxford University, IBM's Almaden Research Center, and Los Alamos. In 2007, a Canadian company D-Wave demonstrated a 16-qubit quantum computer at the Computer History Museum in Mountain View, California. In the following years, several companies announced their successes in developing quantum computers, including Google, IBM, Intel and Rigetti.

As larger and larger quantum computers are built, cryptosystems such as RSA, ElGamal and ECC will be no longer secure, so post-quantum cryptography will be critical for the future of secret communication. In 2016, the National Institute of Standards and Technology (NIST) launched a global project to solicit and select a handful of encryption algorithms with the ability to resist quantum computer attacks. On July 5, 2022, after three rounds of competition and selection, NIST announced four algorithms that will underpin its future cryptography standards. They include one algorithm (CRYSTALS-Kyber) for general encryption and key establishment purposes and three (CRYSTALS-Dilithium, Falcon and Sphincs$+$) for digital signatures. On August 13, 2024, the agency announced three post-quantum cryptography standards: FIPS 203 based on CRYSTALS-Kyber, FIPS 204 based on CRYSTALS-Dilithium, and FIPS 205 based on Sphincs$+$. The fourth standard based on Falcon is on the way.  On November 12, 2024, NIST published a guideline \lq\lq Transition to Post-Quantum Cryptography Standards" which lists detailed route and time table. In fact, many high tech companies and institutions have already completed the transition.

It is well known that CRYSTALS-Kyber and CRYSTALS-Dilithium are based on LWE, Falcon is based on NTRU, and Sphincs$+$ is based on the Hash function. Furthermore, both NTRU and LWE are based on lattices. Lattice-based cryptography was born more or less at the same time of Shor's quantum algorithms for the discrete logarithm problem and the factorization problem. It has been explored as a key candidate for post-quantum cryptography ever since.

This paper will review some key mathematical problems behind the lattice-based cryptography. For a complete picture of how post-quantum cryptography grow up from mathematics we refer to \cite{Hoff08,Zong??}.

\bigskip
\noindent
{\large\bf 2. Computational Complexity Problems}

\smallskip
\noindent
No one can predict the future of the post-quantum cryptography. Up to now, lattice-based cryptosystems play the decisive role. Although different in form, the security of all known lattice-based cryptosystems and algorithms relies on the computational complexity of the following two problems and their variations.

\smallskip
\noindent
{\bf The Shortest Vector Problem (SVP).} {\it Find a shortest nonzero vector in an $n$-dimensional lattice $\Lambda $, i.e., find a nonzero vector
${\bf v}\in \Lambda $ that minimizes the Euclidean norm $\| {\bf v} \|$.}

\smallskip
\noindent
{\bf The Closest Vector Problem (CVP).} {\it Given a vector ${\bf x}\in \mathbb{E}^n$ that is not in $\Lambda $, find a vector ${\bf v}\in \Lambda$ that is closest to ${\bf x}$, i.e., find a vector ${\bf v}\in \Lambda $ that minimizes the Euclidean norm $\| {\bf v}-{\bf x}\|.$}

\smallskip
In fact, the security of Ajtai-Dwork, NTRU and LWE is based on the complexity of SVP and its generalizations, and the security of GGH is based on the complexity of CVP.

\medskip
In 1981, P. van Emde Boas \cite{Boas81} started to study the computational complexity of these problems for classical computers and first proved the following theorem.

\medskip\noindent
{\bf Theorem 2.1.} {\it The closest vector problem is $\mathcal{NP}$-hard.}

\medskip
In that paper he was not able to solve the shortest vector problem. Instead, he made the following conjecture.

\medskip\noindent
{\bf Conjecture 2.1.} {\it The shortest vector problem is $\mathcal{NP}$-hard.}

\medskip
Forty years later, this conjecture is still open today. Meanwhile, research has turned towards randomized reduction and approximation. Unlike deterministic reduction, randomized reduction allows the mapping function to be computable in polynomial time by a probabilistic algorithm. Therefore, the output of the reduction is only required to be correct with sufficiently high probability. In 1997, Ajtai \cite{Ajta98} proved the following theorem.

\medskip\noindent
{\bf Theorem 2.2.} {\it The shortest vector problem is $\mathcal{NP}$-hard under randomized reduction.}

\medskip
In fact, even approximating the shortest lattice vector or the closest lattice vector are not easy. These facts guarantee some important applications. For convenience, let $\ell (\Lambda)$ denote the length of the shortest nonzero vector of $\Lambda$.

\medskip
\noindent
{\bf SVP$_\gamma$.} {\it Let $\Lambda$ be an $n$-dimensional lattice and let $\gamma (n)$ be an approximation factor. Find a nonzero vector}
${\bf v}\in \Lambda $ satisfying $$\| {\bf v} \|\le \gamma (n) \ell (\Lambda ).$$

\medskip
\noindent
{\bf CVP$_\gamma$.} {\it Let $\Lambda$ be an $n$-dimensional lattice, let ${\bf x}$ be a point of $\mathbb{E}^n$, and $\gamma (n)$ be an approximation factor. Find a ${\bf v}\in \Lambda $ satisfying} $$\| {\bf x}, {\bf v} \|\le \gamma (n) \| {\bf x}, \Lambda\|.$$

\medskip
SVP$_\gamma$ and CVP$_\gamma$ have been studied by many authors and a number of results have been achieved. We list two of them here as examples.

\medskip\noindent
{\bf Theorem 2.3 (Khot \cite{Khot05}).} {\it To approximate the shortest vector of an $n$-dimensional lattice within any constant factor $c$ under randomized reduction is $\mathcal{NP}$-hard.}

\medskip\noindent
{\bf Theorem 2.4 (Dinur, Kindler, Raz and Safra \cite{Dinu03}).} {\it To approximate the closest vector of an $n$-dimensional lattice to a given point of $\mathbb{E}^n$ within a factor $n^{c/\log\log n}$, where $c$ is some absolute constant, is $\mathcal{NP}$-hard.}

\medskip
Besides the approximation versions, SVP$_\gamma$ and CVP$_\gamma$ also have the following decision versions.

\medskip
\noindent
{\bf GapSVP$_\gamma$.} {\it Let $\Lambda$ be an $n$-dimensional lattice and let $r$ be a given positive number, the GapSVP$_\gamma$ with approximation factor $\gamma (n)$ asks to decide whether $\ell (\Lambda )\le r$ or $\ell (\Lambda )>  r\!\ \gamma (n)$.}

\medskip
\noindent
{\bf GapCVP$_\gamma$.} {\it Let $\Lambda$ be an $n$-dimensional lattice, let ${\bf x}$ be a point of $\mathbb{E}^n$, and let $r$ be a given positive number. The GapCVP$_\gamma$ with approximation factor $\gamma (n)$ asks to decide whether $\|{\bf v}, {\bf x} \|\le r$ holds for a lattice point ${\bf v}$ or $\|\Lambda, {\bf x}\| >  r\!\ \gamma (n)$.}

\medskip
In 2004, Ajtai introduced a new problem, called the short integer solution (SIS) problem, over random q-ary lattices. He proved that, under certain parameters, solving SIS over a lattice chosen randomly from an easily samplable distribution is at least as hard as approximating the shortest vector problem for any lattice. In fact, the security of both NTRU and LWE does rely on the complexity of approximating versions of the SVP.

\medskip
Since the birth of Shor's quantum algorithms for discrete logarithms and factoring in 1994, in particular since the NIST initiated the post-quantum cryptography competition in 2016, people have tried hard to search for efficient quantum computing algorithms for the shortest vector problem and the closest vector problem (see \cite{Agga25}), or tried to prove that there is no such algorithm. Up to now, none of this effort has succeeded. This failure led to the following conjectures.

\medskip\noindent
{\bf Conjecture 2.2.} {\it There is no polynomial time quantum algorithm which can approximate the shortest vector problem within a polynomial factor.}

\medskip\noindent
{\bf Conjecture 2.3.} {\it There is no polynomial time quantum algorithm which can approximate the closest vector problem within a polynomial factor.}

\medskip
If these two conjectures are correct, the security of the related lattice-based cryptosystems can be guaranteed in the quantum computing era.

\medskip
SVP and CVP have several useful generalizations. We list seven of them as follows. Just like SVP and CVP, one can easily formulate their approximation versions and their decision versions.

\medskip
\noindent
{\bf The Shortest Basis Problem (SBP).} {\it For $B=\{ {\bf a}_1, {\bf a}_2, \ldots,$ ${\bf a}_n\}$ we define
$$\mu (B)=\max_i\| {\bf a}_i\|.$$
To find a basis $B$ of an $n$-dimensional lattice $\Lambda$ such that to minimize the value of $\mu (B)$.}

\medskip
\noindent
{\bf The Covering Radius Problem (CRP).} {\it For a given $n$-dimensional lattice $\Lambda$, to determine the minimum number $r$ such that $rB^n+\Lambda$ is a covering of $\mathbb{E}^n$.}

\medskip
\noindent
{\bf The Quasi Orthogonal Basis Problem (QOBP).} {\it Given an $n$-dimensional lattice $\Lambda$, find a basis $\{{\bf a}_1$, ${\bf a}_2$,
$\ldots $, ${\bf a}_n\}$ to minimize $\prod_{i=1}^n\| {\bf a}_i\|$.}

\medskip
Let $\lambda_i(\Lambda )$ denote the radius of the smallest ball centered at the origin containing $i$ linearly independent lattice vectors. Usually, $\lambda_i(\Lambda )$ is called the $i$-th successive minimum of $\Lambda$.

\medskip
\noindent
{\bf The Successive Minima Problem (SMP).} {\it Given an $n$-dimensional lattice $\Lambda$, find $n$ linearly independent vectors ${\bf v}_1$, ${\bf v}_2$, $\ldots $, ${\bf v}_n$ satisfying $\| {\bf v}_i\|=\lambda_i(\Lambda).$}

\medskip
\noindent
{\bf The Shortest Independent Vectors Problem (SIVP).} {\it Given an $n$-dimensional lattice $\Lambda$, find $n$ linearly independent vectors ${\bf v}_1$, ${\bf v}_2$, $\ldots $, ${\bf v}_n$ satisfying $\| {\bf v}_i\|\le \lambda_n(\Lambda).$}

\medskip
 If $\{ {\bf a}_1, {\bf a}_2, \ldots, {\bf a}_n\}$ be a basis of an $n$-dimensional lattice, let $H_i$ denote the $i$-dimensional hyperplane spaned by the first $i$ vectors of the basis, and let ${\bf b}_i$ denote the projection of ${\bf a}_i$ to $H_{i-1}$. Then we define ${\bf a}_1^*={\bf a}_1$ and ${\bf a}_i^*= {\bf a}_i-{\bf b}_i$ for $i=2, 3, \ldots, n$.

\medskip
\noindent
{\bf The Shortest Diagonal Problem (SDP).} {\it Given an $n$-dimen- sional lattice $\Lambda$, find a basis $\{{\bf a}_1$, ${\bf a}_2$,
$\ldots $, ${\bf a}_n\}$ to minimize $\sum_{i=1}^n\|{\bf a}_i^*\|^2.$}

\medskip
\noindent
{\bf The Densest Sublattice Problem (DSP).} {\it Let $k$ be an integer satisfying $k<n$. Given an $n$-dimensional lattice, find $k$ linearly independent lattice vectors ${\bf b}_1$, ${\bf b}_2$, $\ldots$, ${\bf b}_k$ to minimize the $k$-dimensional volume of $\{\lambda_1{\bf b}_1+\lambda_2{\bf b}_2+\ldots +\lambda_k{\bf b}_k:\ 0\le \lambda_i\le 1\}$.}

\medskip
These problems have been studied extensively for classical computers (see \cite{Micc02,Wang23,Zhang20}). For quantum computers, little is known.

\bigskip
\noindent
{\large\bf 3. Geometric Problems}

\smallskip
\noindent
Let $B^n$ denote the $n$-dimensional unit ball $\{ {\bf x}: \sum x_i^2\le 1\}$ in $\mathbb{E}^n$ with volume $\omega_n$ and let $\Lambda$ denote an $n$-dimensional lattice with determinant $d(\Lambda )$. We call $B^n+\Lambda=\{ B^n+{\bf v}_i: {\bf v}_i\in \Lambda\}$ a ball packing\footnote{In discrete geometry, it is called sphere packing rather than ball packing.} if the interiors of the balls are disjoint. Let $\delta^*(B^n)$ denote the density of the densest lattice ball packings. In other words,
$$\delta^*(B^n)=\max_{\Lambda} {{\omega_n}\over {d(\Lambda )}},$$
where the maximum is over all the lattices $\Lambda $ such that $B^n+\Lambda$ are ball packings.

Let $\ell (\Lambda )$ denote the length of the shortest nonzero vectors of $\Lambda$ and take $r={1\over 2}\ell (\Lambda )$. It is easy to see that $rB^n+\Lambda $ is a lattice ball packing in $\mathbb{E}^n$. Then, the shortest vector problem can be reformulated in terms of ball packing.

\medskip
\noindent
{\bf SVP in Ball Packing.} {\it For a given $n$-dimensional lattice $\Lambda $, find the largest number $r$ such that $rB^n+\Lambda$ is a ball packing and the corresponding balls that touch $rB^n$ at its boundary.}

\medskip
In fact, based on the previous discussion, one can deduce the following connection between the length $\ell (\Lambda )$ of the shortest nonzero vector of a lattice $\Lambda $ and the lattice ball packing densities $\delta^*(B^n)$.

\medskip
\noindent
{\bf Theorem 3.1.} {\it Let $\Lambda$ be an $n$-dimensional lattice and let $\omega_n$ denote the volume of $B^n$. We have}
$$\ell (\Lambda)\le 2 \sqrt[n]{{\rm det}(\Lambda )\cdot\delta^*(B^n)/\omega_n}.$$

\medskip
Therefore, the following problem is fundamental for the SVP.

\medskip
\noindent
{\bf Problem 3.1.} {\it To determine or estimate the values of $\delta^*(B^n).$}

\medskip
This problem and its generalizations have been studied by many prominent mathematicians, including Kepler, Gauss, Hermite, Voronoi and Minkowski. However, our knowledge in this field is still very limited. Some key results about $\delta^*(B^n)$ are summarized in the following table.

{\large
$$\begin{tabular}{|c|c|c|c|c|}
\hline
{\small n}&{\small 2}&{\small 3}&{\small 4}&{\small 5}\\
\hline
{\small $\delta^*(B^n)$}& ${{\pi }\over {\sqrt{12}}}$& ${{\pi }\over {\sqrt{18}}}$& ${{\pi^2}\over {16}}$& ${{\pi^2}\over {15\sqrt{2}}}$\\
\hline
\mbox{${{\rm Author}\atop {\rm Date}}$} & \mbox{${{\rm Lagrange}\atop {\rm 1773}}$} & \mbox{${{\rm Gauss}\atop {\rm 1831}}$} & \mbox{${{{\rm Korkin,}\atop {\rm Zolotarev}}\atop {\rm 1872}}$} & \mbox{${{{\rm Korkin,}\atop {\rm Zolotarev}}\atop {\rm 1877}}$} \\
\hline
\hline
{\small n}&{\small 6}&{\small 7}&{\small 8}&{\small 24}\\
\hline
{\small $\delta^*(B^n)$}& ${{\pi^3}\over {48\sqrt{3}}}$& ${{\pi^3}\over {105}}$& ${{\pi^4}\over {384}}$& ${{\pi^{12}}\over {12!}}$\\
\hline
\mbox{${{\rm Author}\atop {\rm Date}}$} & \mbox{${{\rm Blichfeldt}\atop {\rm 1925}}$} & \mbox{${{\rm Blichfeldt}\atop {\rm 1926}}$} & \mbox{${{\rm Blichfeldt}\atop {\rm 1934}}$} & \mbox{${{\rm Cohn,\ Kumar}\atop {\rm 2009}}$}\\
\hline

\end{tabular}$$}

In general dimensions, we have
$$cn^22^{-n}\le \delta^*(B^n)\le 2^{-0.599n(1+o(1))}$$
for a suitable positive constant $c$, where a weaker lower bound was first proved by Minkowski in 1905, then improved and generalized by E. Hlawka, C. L. Siegel, H. Davenport, C. A. Rogers, W. M. Schmidt, B. Klartag and others (see \cite{Klar25}), and the upper bound was proved by G. A. Kabatjanski and V. I. Leven$\check{\rm s}$tein in 1978 (see \cite{Zong99}).

There are hundreds of papers on ball packings, employing methods and tools from various fields of mathematics. For example, the surprising works of M. Viazovska and her coauthors \cite{Cohn17,Viaz17} in $\mathbb{E}^8$ and $\mathbb{E}^{24}$, respectively. As well, many fascinating problems on ball packings are still open. Here we list one of them as an example.

\medskip\noindent
{\bf Problem 3.2.} {\it Determine the asymptotic order of $\delta^*(B^n)$, if it does exist.}

\medskip
Assume that $\Lambda$ is an $n$-dimensional lattice in $\mathbb{E}^n$. For every point ${\bf x}\in \mathbb{E}^n$, we define the distance between ${\bf x}$ and its closest lattice point ${\bf v}\in \Lambda$ as $d ({\bf x},\Lambda)$. Then, we define
$$\rho (\Lambda )=\max_{{\bf x}\in \mathbb{E}^n} d({\bf x}, \Lambda).$$
It is easy to see that $\rho(\Lambda)B^n+\Lambda$ is a covering of $\mathbb{E}^n$. In fact, $\rho (\Lambda)$ is the smallest radius $\rho$ such that $\rho B^n+\Lambda$ is a covering of $\mathbb{E}^n$.

\medskip
\noindent
{\bf CVP in Ball Covering.} {\it Given an $n$-dimensional lattice $\Lambda$, find the smallest number $\rho$ such that $\rho B^n+\Lambda$ is a covering of $\mathbb{E}^n$. For any ${\bf x}\in \mathbb{E}^n$, find a lattice point ${\bf v}\in \rho B^n+{\bf x}$.}

\medskip
This version is closely related to the CRP and the SDP introduced in the previous section. Clearly, finding a lattice point ${\bf v}\in \rho B^n+{\bf x}$ is slightly simpler than the CVP. However, this covering model can illustrate the fundamental difficulty of the CVP. Simple example can show that there is no upper bound for $\rho (\Lambda )$ in terms of ${\rm det}(\Lambda)$ and $n$.

Let $\theta^* (B^n)$ denote the density of the thinnest lattice ball covering of $\mathbb{E}^n$. In other words,
$$\theta^*(B^n)=\min_{\Lambda} {{\omega_n}\over {d(\Lambda )}},$$
where the minimum is over all the lattices $\Lambda $ such that $B^n+\Lambda$ are ball coverings of $\mathbb{E}^n$.
As a counterpart to Theorem 3.1, we have the following relation between $\rho (\Lambda )$ and $\theta^*(B^n)$.

\medskip
\noindent
{\bf Theorem 3.2.} {\it Let $\Lambda$ be an $n$-dimensional lattice and let $\omega_n$ denote the volume of $B^n$. We have}
$$\rho (\Lambda)\ge \sqrt[n]{{\rm det}(\Lambda )\cdot\theta^*(B^n)/\omega_n}.$$

\medskip
Therefore, the following problem is fundamental for CVP, CRP and SDP.

\medskip\noindent
{\bf Problem 3.3.} {\it Determine the asymptotic order of $\theta^*(B^n)$, if it does exist.}

\medskip
Ball covering, in certain sense, is regarded as a dual concept of ball packing. In fact, they are not much related. Up to now, the known exact results about $\theta^* (B^n)$ can be summarized as follows.

{\large
$$\begin{tabular}{|c|c|c|c|c|}
\hline
{\small n}& {\small 2}& {\small 3}&{\small 4}&{\small 5}\\
\hline
{\small $\theta^*(B^n)$}& ${{2\pi }\over {3\sqrt{3}}}$& ${{5\sqrt{5}\pi }\over {24}}$& ${{2\pi^2}\over {5\sqrt{5}}}$& ${{245\sqrt{35}\pi^2}\over {3888\sqrt{3}}}$\\
\hline
\mbox{${{\rm Author}\atop {\rm Date}}$}& \mbox{${{\rm Kersshner}\atop {\rm 1939}}$}& \mbox{${{\rm Bambah}\atop {\rm 1954}}$}& \mbox{${{\rm Delone,\ Ryskov}\atop {\rm 1963}}$}& \mbox{${{{\rm Ryskov,}\atop {\rm Baranovskii}}\atop {\rm 1975}}$}\\
\hline
\end{tabular}$$}

\medskip
In general dimensions, there are three constants $c_1$, $c_2$ and $c_3$ such that
$$c_1 n\le \theta^*(B^n)\le c_2n(\log_en)^{c_3},$$
where the lower bound was achieved by H. S. M. Coxeter, L. Few and Rogers in 1959, and the upper bound was discovered by Rogers in 1959 (see \cite{Roge64}).

\medskip
One may realize that there are very few concrete results on ball covering in the past half a century, particularly compared to ball packing. It is fascinating to notice that, unlike the packing case, the thinnest lattice ball covering in $\mathbb{E}^8$ is not achieved by the $E_8$ lattice. At least, the $A_8^*$ lattice provides a ball covering with a density thinner than the $E_8$ lattice. Therefore, the following problem is fascinating for mathematicians.

\medskip\noindent
{\bf Problem 3.4.} {\it Determine the values of $\theta^*(B^8)$ and $\theta^*(B^{24})$.}

\medskip
Let $\mathcal{L}_n$ denote the family of all $n$-dimensional lattices. In 1950, Rogers defined and studied
$$\phi^*(B^n)=\min_{\Lambda\in \mathcal{L}_n}{{2\rho (\Lambda)}\over {\ell (\Lambda)}},$$
where $\ell (\Lambda)$ is the length of the shortest nonzero vectors of $\Lambda$ and $\rho (\Lambda)$ is the maximum distance between a point ${\bf x}\in \mathbb{E}^n$ and its closest lattice point. In other words, $\phi^*(B^n)$ is the optimal ratio of the covering radius to the packing radius of an $n$-dimensional lattice. Therefore, it is a bridge between SVP and CVP.

From the intuitive point of view, one may think that $\phi^*(B^n)$ can be arbitrarily large when $n\rightarrow \infty$. Surprisingly, Rogers proved by a reduction method that
$$\phi^*(B^n)\le 3$$
holds in every dimension. In 1972, via mean value techniques developed by Rogers and Siegel, G. L. Butler improved Rogers' upper bound to
$$\phi^*(B^n)\leq 2+o(1).$$

\medskip\noindent
{\bf Problem 3.5.} {\it Determine the value of $$\lim_{n\to\infty }\phi^*(B^n),$$ if it does exist.}

\medskip
It follows from Rogers' upper bound that, for a proportion of $n$-dimensional lattices, the longest distance in CVP is only a constant multiple of the length of the SVP. Recent years, this idea has been applied to cryptographic analysis by Micciancio and others.

In the 1980s, several mathematicians studied $\phi^*(B^n)$ from different respects. Up to now, we have the following exact results.

$$\begin{tabular}{|c|c|c|c|c|c|}
\hline
$n$ &{\small $2$} &{\small $3$} &{\small $4$} &{\small $5$} \\
\hline
{\small $\phi^*(B^n)$} &{\small $2/\sqrt3$ }&{\small $\sqrt{5/3}$} &{\small $\sqrt{2\sqrt3}(\sqrt3-1)$} &{\small $\sqrt{{3\over 2}+{\sqrt{13}\over 6}}$}\\
\hline
\mbox{${{\rm Author}\atop {\rm Date}}$} &  & \mbox{${{\rm Boroczky}\atop {\rm 1986}}$} & \mbox{${{\rm Horvath}\atop {\rm 1982}}$} & \mbox{${{\rm Horvath}\atop {\rm 1986}}$}\\
\hline
\end{tabular}$$

\medskip
For mathematicians, the following two problems are particular interesting.

\medskip\noindent
{\bf Problem 3.6.} {\it Determine the values of $\phi^*(B^8)$ and $\phi^*(B^{24})$.}

\medskip\noindent
{\bf Problem 3.7.} {\it Is there a dimension $n$ such that $\phi^*(B^n)\ge 2 ?$}

\medskip
If one can improve Butler's upper bound to $\phi^*(B^n)\le 2-c$, where $c$ is a positive constant, the lower bound for $\delta^*(B^n)$ will be improved to
$$\delta^*(B^n)\ge (2-c)^{-n}.$$
If a dimension $n$ can be found such that $\phi^*(B^n)\ge 2$, then we have
$$\delta^*(B^n)\not=\delta (B^n),$$
where $\delta (B^n)$ is the density of the densest ball packing in $\mathbb{E}^n$. It is easy to see that $\phi^*(B^n)$ can be generalized from the ball to arbitrary centrally symmetric convex bodies. For more on $\phi^*(B^n)$ and its generalizations, we refer to Zong \cite{Zong02}.

\medskip
There is another important notion which is closely related to both the shortest vector problem and the closest vector problem, the Dirichlet-Voronoi cell of $\Lambda$:
$$D=\left\{ {\bf x}:\ {\bf x}\in \mathbb{E}^n,\ \mbox{$\langle {\bf x}, {\bf v}\rangle \le {1\over 2}\langle {\bf v}, {\bf v}\rangle$ for all ${\bf v}\in \Lambda\setminus \{ {\bf o}\}$}\right\}.$$
Roughly speaking, $D$ is the set of points that closer to the origin than any other lattice point. Clearly, $D$ is a centrally symmetric polytope such that $D+\Lambda$ is a tiling of $\mathbb{E}^n$. Therefore, the closest vector problem can be reformulated as:

\medskip
\noindent
{\bf CVP in D-V Cell.} {\it Let $\Lambda$ be an $n$-dimensional lattice and ${\bf x}$ be an arbitrary point of $\mathbb{E}^n$. If $D$ is the Dirichlet-Voronoi cell of $\Lambda$, find a lattice point ${\bf v}$ satisfying ${\bf x}\in D+{\bf v}$.}

\medskip
Clearly, the D-V cell encodes a lot of information of the lattice. Let $\Gamma $ denote the boundary of $D$, then one can easily deduce that
$$\ell (\Lambda )=2 \min_{{\bf x}\in \Gamma }\| {\bf o}, {\bf x} \|$$
and
$$\rho (\Lambda)=\max_{{\bf x}\in \Gamma} \| {\bf o}, {\bf x} \|.$$
Therefore, to understand the D-V cells of lattices will be helpful to solve SVP, CVP and CRP.

We end this section with two well known problems about the Dirichlet-Voronoi cells of lattices. As usually, a parallelotope is a convex polytope which can form a lattice tiling of $\mathbb{E}^n$.

\medskip\noindent
{\bf Problem 3.8.} {\it When $n\ge 6$, classify all Dirichlet-Voronoi cells of the $n$-dimensional lattices, i.e., determine their geometric shapes.}

\medskip\noindent
{\bf Voronoi's Conjecture.} {\it Every parallelotope is a linear image of a lattice Dirichlet-Voronoi cell.}

\medskip
When $n\le 5$, both Problem 3.8 and Voronoi's conjecture have been solved. The Dirichlet-Voronoi cell has been applied to lattice-based cryptography by Micciancio and others since 2010.

\bigskip
\noindent
{\large\bf 4. Arithmetic Problems}

\smallskip
\noindent
Let $\Lambda $ be a lattice with a basis $\{ {\bf a}_1,$ ${\bf a}_2,$ $\ldots ,$ ${\bf a}_n\}$, where ${\bf a}_i=(a_{i1}, a_{i2},$ $\ldots ,$ $a_{in})$, and let $A$ denote the $n\times n$ matrix with entries $a_{ij}$. Then, the lattice can be expressed as
$$\Lambda =\left\{{\bf z}A:\ {\bf z}\in \mathbb{Z}^n\right\}$$
and the norms of the lattice vectors  can be expressed as a positive definite quadratic form
$$Q({\bf z})=\langle {\bf z}A, {\bf z}A\rangle ={\bf z}AA'{\bf z}',$$
where $A'$ and ${\bf z}'$ indicate the transposes of $A$ and ${\bf z}$, respectively. Assume that
$$Q({\bf x})=\sum_{1\le i,j\le n}c_{ij}x_ix_j={\bf x}C{\bf x}'$$
is a positive definite quadratic form of $n$ variables, where $c_{ij}=c_{ji}$ and $C$ is the symmetric matrix with entries $c_{ij}$. It is known that there is an $n\times n$ matrix $A$ satisfying $C=AA'$. Then the quadratic form also produces a lattice
$$\Lambda =\left\{{\bf z}A:\ {\bf z}\in \mathbb{Z}^n\right\}.$$
Therefore, there is a nice correspondence between lattices and positive definite quadratic forms. For convenience, we define
$$m(Q)=\min_{{\bf z}\in \mathbb{Z}^n\setminus \{ {\bf o}\}} Q({\bf z}).$$
Then, the shortest vector problem and its approximation version SVP$_\gamma$ can be reformulated as follows.

\medskip
\noindent
{\bf SVP in Quadratic Form.} {\it Find a nonzero vector ${\bf z}\in \mathbb{Z}^n$ that minimizes the positive definite quadratic form $Q({\bf z})$. In other words, to find an integer solution of}
$$Q({\bf z})= m(Q).$$

\medskip
\noindent
{\bf SVP$_\gamma$ in Quadratic Form.} {\it Let $\gamma\ge 1$ be a positive constant or a polynomial of $n$. Find a nonzero vector ${\bf z}\in \mathbb{Z}^n$ satisfying}
$$Q({\bf z})\le \gamma^2m(Q).$$

\medskip
Let ${\rm dis} (Q)$ be the discriminant of the quadratic form $Q({\bf x})$ and let $\mathcal{Q}_n$ denote the family of all positive definite quadratic forms in $n$ variables. Hermite defined and studied
$$\gamma_n=\sup_{Q\in \mathcal{Q}_n} {{m(Q)}\over {\sqrt[n]{ {\rm dis}(Q)}}}.$$
Usually, $\gamma_n$ is called Hermite's constant. Similar to Theorem 3.1, we have

\medskip\noindent
{\bf Theorem 4.1.} {\it For any positive definite quadratic form $Q({\bf x})$ of $n$ variables, we have}
$$m(Q)\le \gamma_n {\sqrt[n]{ {\rm dis}(Q)}}.$$

\medskip
Therefore, the following problem is basic to understand SVP and SVP$_\gamma$.

\medskip\noindent
{\bf Problem 4.1.} {\it Determine or estimate Hermite's constants $\gamma_n$.}

\medskip
These constants are closely related to the densities $\delta^*(B^n)$ of the densest lattice ball packings. Since $\ell (\Lambda )=\sqrt{m(Q)}$ and ${\rm dis}(Q)={\rm det}(\Lambda )^2$, one can easily deduce that
$$\delta^*(B^n)={{\omega_n\gamma_n^{n/2}}\over {2^n}},$$
where $\omega_n$ is the volume of the $n$-dimensional unit ball $B^n$. In fact, all the known exact results about $\delta^*(B^n)$ (except $\delta^*(B^{24})$) were derived from $\gamma_n$.

{\large
$$\begin{tabular}{|c|c|c|c|c|}
\hline
{\small n}&{\small 2}&{\small 3}&{\small 4}&{\small 5}\\
\hline
$\gamma_n$& {\small $2/\sqrt3$ }&{\small $\sqrt[3]{2}$}& {\small $\sqrt{2}$} & {\small $\sqrt[5]{8}$}\\
\hline
\mbox{${{\rm Author}\atop {\rm Date}}$} & \mbox{${{\rm Lagrange}\atop {\rm 1773}}$} & \mbox{${{\rm Gauss}\atop {\rm 1831}}$} & \mbox{${{{\rm Korkin,}\atop {\rm Zolotarev}}\atop {\rm 1872}}$} & \mbox{${{{\rm Korkin,}\atop {\rm Zolotarev}}\atop {\rm 1877}}$} \\
\hline
\hline
{\small n}&{\small 6}&{\small 7}&{\small 8}&{\small 24}\\
\hline
$\gamma_n$&{\small $\sqrt[6]{64\over 3}$}& {\small $\sqrt[7]{64}$} & {\small $2$}& {\small $4$}\\
\hline
\mbox{${{\rm Author}\atop {\rm Date}}$} &\mbox{${{\rm Blichfeldt}\atop {\rm 1925}}$} & \mbox{${{\rm Blichfeldt}\atop {\rm 1926}}$} & \mbox{${{\rm Blichfeldt}\atop {\rm 1934}}$} & \mbox{${{\rm Cohn,\ Kumar}\atop {\rm 2009}}$}\\
\hline
\end{tabular}$$}

\medskip
In 1953, R. A. Rankin \cite{Rank53} introduced a generalization of Hermite's constant. Let $k$ be an integer, $1\le k\le n-1$, and let $m_k(Q)$ denote the lower bound of any principal minor of order $k$ of any form equivalent to $Q({\bf x})$. He defined and studied
$$\gamma_{n,k}=\sup_{Q\in \mathcal{Q}_n} {{m_k(Q)}\over { {\rm dis}(Q)^{k/n}}}.$$
Clearly, we have $\gamma_{n,1}=\gamma_n$. Usually, $\gamma_{n,k}$ are called Rankin's constants.

\medskip\noindent
{\bf Problem 4.2.} {\it Determine or estimate Rankin's constants $\gamma_{n,k}$.}

\medskip
Up to now, our knowledge about these constants is very limited. In fact, the only known nontrivial exact result is $\gamma_{4,2}=3/2$, which was discovered by Rankin himself. Twenty years ago, Rankin's constant leaded P. Nguyen and others to introduce the densest sublattice problem (see the DSP at the end of Section 2), a generalization of the shortest vector problem. This new problem has been studied by Micciancio, Nguyen and others. It has important applications to blockwise lattice reduction generalizing LLL and Schnorr's algorithm.

\medskip
Assume that $\Lambda =\{{\bf z}A:\ {\bf z}\in \mathbb{E}^n\}$ is an $n$-dimensional lattice in $\mathbb{E}^n$, where $A$ is a nonsingular $n\times n$ matrix. For any point ${\bf p}={\bf y}A\in \mathbb{E}^n$ and ${\bf v}={\bf z}A\in \Lambda$, we have
$$\| {\bf p}-{\bf v}\|=\| ({\bf y}-{\bf z})A\|=\sqrt{Q({\bf y}-{\bf z})}.$$
Therefore, the closest vector problem is equivalent to the following problem.

\medskip
\noindent
{\bf CVP in Quadratic Form.} {\it Given a positive definite quadratic form $Q({\bf x})$ and a vector ${\bf y}$, find an integer vector ${\bf z}\in \mathbb{Z}^n$ that minimizes $Q({\bf y}-{\bf z})$.}

\medskip
Let $C$ denote the unit cube $\{(x_1, x_2, \ldots , x_n):\ 0\le x_i< 1\}$, let $\Lambda $ be the lattice corresponding to $Q({\bf x})$, and define
$$\rho (Q) =\sqrt{\max_{{\bf y}\in C}\min_{{\bf z}\in \mathbb{Z}^n}Q({\bf y}-{\bf z})}.$$
It can be verified that $\rho (Q)$ is the smallest number $\rho$ such that $\rho B^n+\Lambda $ is a ball covering of $\mathbb{E}^n$. Therefore, the CRP can be reformulated as follows.

\medskip
\noindent
{\bf CRP in Quadratic Form.} {\it Given a positive definite quadratic form $Q({\bf x})$, determine the value of}
$$\max_{{\bf y}\in C}\min_{{\bf z}\in \mathbb{Z}^n}Q({\bf y}-{\bf z}).$$

\medskip
Clearly, we have
$$\theta^*(B^n)=\min_{Q\in \mathcal{Q}_n}{{\omega_n\rho (Q)^n}\over {\sqrt{{\rm dis}(Q)}}}.$$
In fact, most of the known exact results about $\theta^*(B^n)$ were also achieved by studying quadratic forms.

\medskip
If one knows that $\Lambda$ has an orthogonal basis $\{{\bf a}_1$, ${\bf a}_2$, $\ldots $, ${\bf a}_n\}$, then both SVP and CVP are trivial. Of course, most lattices have no orthogonal bases. Nevertheless, every lattice has some relatively good bases in certain sense. Correspondingly, every positive definite quadratic form has a comparatively good equivalent form. This is the philosophy of reduction theory. In the history, reduction theory was first developed for quadratic forms rather than for lattices.

In 1773, Lagrange proved that every positive definite binary quadratic form $Q({\bf x})={\bf x}C{\bf x}'$ is equivalent to one satisfying
$$\left\{
\begin{array}{ll}
c_{11}\le c_{22},& \mbox{}\\
0\le 2c_{12}\le c_{11},& \mbox{}
\end{array}\right.$$
which marked the birth of the reduction theory. In other words, every two-dimensional lattice has a basis $\{{\bf a}_1, {\bf a}_2\}$ such that the angle between ${\bf a}_1$ and ${\bf a}_2$ is at least $\pi/3$ and at most $\pi/2$. Then, one can deduce that $\gamma_2=2/\sqrt{3}$ and $\delta^*(B^2)=\pi/\sqrt{12}.$

As a generalization of Lagrange's pioneering work, in 1905 Minkowski discovered the following reduction: For convenience, we denote the greatest common divisor
of $k$ integers $z_1$, $z_2$, $\ldots ,$ $z_k$ by $[z_1,z_2, \ldots, z_k]$. A positive definite
quadratic form $Q({\bf x})= {\bf x}C{\bf x}'$ is said to be Minkowski reduced, if
$$c_{1j}\ge 0, \quad j=2,\ 3,\ \ldots ,\ n,$$
and
$$Q({\bf z})\ge c_{ii}, \quad i=1,\ 2,\ \ldots ,\ n,$$
for all integer vectors ${\bf z}=(z_1, z_2, \ldots , z_n)$ satisfying $[z_i, z_{i+1}, \ldots,$ $z_n]=1.$

\medskip
Then, he proved the following theorem.

\medskip
\noindent
{\bf Theorem 4.2.} {\it Every positive definite quadratic form is equivalent to a Minkowski reduced one.}

\medskip
It is easy to see that the first basis vector in the corresponding lattice of a Minkowski reduced form is the shortest nonzero lattice vector. However, to deduce a Minkowski reduced basis is not easy. Reduction theory has been developed by Seeber, Gauss, Hermite, Korkin, Zolotarev, Minkowski, Voronoi and many modern authors (see \cite{Nguy09, Zong99}). Several different reductions such as the Korkin-Zolotarev reduction and the Lenstra-Lenstra-Lovazs reduction have been discovered. Clearly, reduction methods are the key tools for the security analysis of lattice-based cryptography.

\medskip\noindent
{\bf Problem 4.3.} {\it To develop efficient reductions for quantum computing or to prove their non-existence.}

\medskip
For classical computer, the Lestra-Lenstra-Lovazs reduction was very successful to approximate the shortest lattice vector. For quantum computing our knowledge is too limited.

\medskip
Let $Q({\bf x})={\bf x}C{\bf x}'$ be a positive definite quadratic form of $n$ variables. We define
$$\omega (Q)=\sqrt{{{dis({Q})}\over {c_{11}c_{22}\cdots c_{nn}}}}$$
and
$$\omega_n=\min_{Q\in \mathbb{Q}_n}\max_{Q'}\omega (Q),$$
where the inner maximum is over all forms which are equivalent to $Q({\bf x})$. Let $B=\{ {\bf a}_1, {\bf a}_2, \ldots , {\bf a}_n\}$ be a basis of $\Lambda$. We define
$$\varpi (B)={{d(\Lambda)}\over {\|{\bf a}_1\| \|{\bf a}_2\|\cdots \|{\bf a}_n\|}}$$
and
$$\varpi_n=\min_{\Lambda\in \mathbb{L}_n}\max_{B'}\varpi (B'),$$
where the inner maximum is over all bases of $\Lambda $. Usually, $\varpi (B)$ is called the orthogonality defect of the basis. It is easy to see that
$$\omega_n=\varpi_n.$$
Although $\varpi_n$ is more natural than $\omega_n$, in the history $\omega_n$ has been studied more frequently.

Clearly, to determine the values of $\omega (Q)$ or $\varpi (B)$ are equivalent to the quasi orthogonal basis problem introduced at the end of Section 2. Therefore, the following problem is basic for lattice-based cryptography.

\medskip\noindent
{\bf Problem 4.4.} {\it To determine or estimate the values of $\omega_n$.}

\medskip
In the history, when $n\le 5$ these numbers have been studied by B. L. van der Waerden, E. S. Barnes and others for Minkowski reduced bases. For the general case, our knowledge is very limited.

\medskip\noindent
{\bf Concluding Remarks.} To design good lattice-based cryptosystems or to analysis their security, deep understanding of these mathematical problems is inevitable. Therefore, mathematicians can play important roles in post-quantum cryptography by making contributions to these fundamental problems.

The author acknowledges that the background materials and the references of this paper are similar to those of \cite{Zong??}. While \cite{Zong??} tried to show the general scientists that mathematics is fundamental for post-quantum cryptography, this article try to tell mathematicians what they may contribute to this urgent technology.

\medskip\noindent
{\bf Acknowledgement.} For their constant supports, the author is grateful to G. Tian and X. Wang. This work is supported by the National Natural Science Foundation of China (NSFC12226006, NSFC11921001) and the Natural Key Research and Development Program of China (2018YFA0704701).

\bigskip
\noindent
Chuanming Zong, Center of Applied Mathematics, Tianjin University, Tianjin 300072, China.

\smallskip
\noindent
Email: cmzong@math.pku.edu.cn.

\end{document}